\documentclass{amsart}

\usepackage{amsmath,amssymb,amsthm,latexsym}

\newtheorem{theorem}{Theorem}
\newcommand{\bt}{\begin{theorem}}
\newcommand{\et}{\end{theorem}}
\newtheorem{lemma}{Lemma}
\newcommand{\bl}{\begin{lemma}}
\newcommand{\el}{\end{lemma}}
\newtheorem{corollary}{Corollary}
\newcommand{\bc}{\begin{corollary}}
\newcommand{\ec}{\end{corollary}}
\newcommand{\bconj}{\begin{conjecture}}
\newcommand{\econj}{\end{conjecture}}
\newtheorem{problem}{Problem}
\newcommand{\bprob}{\begin{problem}}
\newcommand{\eprob}{\end{problem}}
\newcommand{\beq}{\begin{equation}}
\newcommand{\eeq}{\end{equation}}
\newcommand{\benum}{\begin{enumerate}}
\newcommand{\eenum}{\end{enumerate}}
\newcommand{\N}{\ensuremath{ \mathbf N }}
\newcommand{\Z}{\ensuremath{\mathbf Z}}

\newcommand{\R}{\ensuremath{\mathbf R}}

\newcommand{\mcr}{\ensuremath{ \mathcal R}}

\newcommand{\mbn}{\ensuremath{ \mathbf n}}

\newcommand{\mbx}{\ensuremath{ \mathbf x}}
\newcommand{\mby}{\ensuremath{ \mathbf y}}
\newcommand{\mbz}{\ensuremath{ \mathbf z}}

\DeclareMathOperator{\card}{\text{card}}

\newcommand{\bmat}{\left(\begin{matrix}}
\newcommand{\emat}{\end{matrix}\right)}
\newcommand{\bsmallmat}{\left(\begin{smallmatrix}}
\newcommand{\esmallmat}{\end{smallmatrix}\right)}

\title{Intersections of sumsets in additive number theory}

\author{Melvyn B.  Nathanson}
\address{Department of Mathematics\\Lehman College (CUNY)\\Bronx, NY 10468}
\email{melvyn.nathanson@lehman.cuny.edu}

\date{\today}

\subjclass[2000]{11B13, 11B05, 11B75,  11P70, 22D99}
\keywords{Sumsets in semigroups, Haar measure of sumsets, additive number theory, 
combinatorial number theory}
\thanks{Supported in part by  PSC-CUNY Research Award Program grant 66197-00 54.}

\begin{document}

\begin{abstract}
Let $A$ be a subset of  an additive abelian semigroup $S$ and let $hA$ be the $h$-fold 
sumset of $A$.   The following question is considered:
Let $(A_q)_{q=1}^{\infty}$ be a strictly decreasing sequence of sets in 
$S$ and let  $A = \bigcap_{q=1}^{\infty} A_q$.  
When does one have
\[
hA = \bigcap_{q=1}^{\infty} hA_q
\]
for some or all $h \geq 2$?
\end{abstract}

\date{\today}

\maketitle

\section{The sumset intersection problem} 
Let $\N$, $\N_0$, and \Z\ denote, respectively, the sets of positive integers, nonnegative integers, 
and integers. 

 Let $S$ be an additive abelian semigroup and let $A$ be a subset of $S$.  
For every positive integer $h$, the $h$-fold sumset of $A$ is the set of all sums of $h$ not 
necessarily distinct elements of $A$, that is, 
\[
hA = \{a_1 + \cdots + a_h: a_i \in A \text{ for all } i = 1,\ldots, h\}.  
\]  
For all $x \in S$, we denote by $r_{A,h}(x)$ the number of representations of $x$ as a sum of 
$h$ elements of $A$, that is, 
\[
r_{A,h}(x) = \card\left\{ (a_1,\ldots, a_h) \in A^h : a_1 + \cdots + a_h = x\right\}.
\]

The subset $A$ is a \emph{basis of order $h$} for $S$ if $hA=S$, or,  
equivalently, if $r_{A,h}(x) \geq 1$ for all $x \in S$.    
The subset $A$ is  a \emph{basis} for $S$ if  $A$ is a  basis of order $h$ 
for $S$ for some positive integer $h$.  
The subset $A$ is a \emph{nonbasis of order $h$} for $S$ 
if  $A$ is not a basis of order $h$ for $S$, that is, if $hA\neq S$.  
The set $A$  is a \emph{nonbasis} for $S$ if $hA\neq S$ for every positive integer $h$.  

The sequence of sets $(A_q)_{q=1}^{\infty}$ is \emph{decreasing}\index{decreasing sequence} 
if $A_q \supseteq A_{q+1}$ for all $q \in \N$.  The decreasing sequence $(A_q)_{q=1}^{\infty}$  is
\emph{strictly decreasing}\index{decreasing sequence!strictly} 
if $ A_q \neq A_{q+1}$ for all $q \in \N$ and  
\emph{asymptotically strictly decreasing}\index{sequence!asymptotically strictly decreasing} 
if $ A_q \neq A_{q+1}$ for infinitely many $q$. 
A decreasing sequence is asymptotically strictly decreasing if and only if it is not eventually constant.  
Every asymptotically strictly decreasing sequence contains a strictly decreasing subsequence. 

Let $(A_q)_{q=1}^{\infty}$ be an asymptotically strictly decreasing sequence of sets 
in the semigroup $S$ and let 
\[
A = \bigcap_{q=1}^{\infty} A_q.
\]
For all $q \in \N$, 
we have $A \subseteq A_q$ and so $hA \subseteq hA_q$.  It follows that  
\beq           \label{intersect:h-intersect}
hA =  h\left( \bigcap_{q=1}^{\infty} A_q \right) \subseteq \bigcap_{q=1}^{\infty} hA_q.
\eeq

This inclusion can be proper.  For example, in the semigroup \Z\ of integers, 
for every positive integer $q$, we define 
\beq           \label{intersect:Asharp}
 A^{\sharp}_q =  \{r\in \Z: |r| \geq q \}. 
\eeq
Let $h \geq 2$.   For every nonnegative integer $n$, we have  
$\pm(n+q) \in A^{\sharp}_q $ and $\pm(h-1)q\in A^{\sharp}_q $.  
It follows that  
 \[
n =  (n+q) + (h-2)q + (-(h-1)q) \in h A^{\sharp}_q   
 \]
and 
\[
-n = -(n+q) + (h-2)(-q) + ((h-1)q)  \in h A^{\sharp}_q   
\]
and so $h A^{\sharp}_q = \Z$.  

Let $A$ be any finite set of integers.   
The sumset $hA$ is finite for all $h \in \N$ and so $hA \neq \Z$.  
The asymptotically strictly decreasing sequence of sets $(A_q)_{q=1}^{\infty}$ defined by 
\[
A_q = A \cup A^{\sharp}_q = A \cup \{r\in \Z: |r| \geq q \} 
\]
satisfies 
$ A = \bigcap_{q=1}^{\infty} A_q$.  
Because $hA_q=\Z$ for all $q \in \N$ , we have 
\[
 hA  \neq \Z =  \bigcap_{q=1}^{\infty} hA_q.  
\]

\bprob 
We want to know when equality replaces inclusion in relation~\eqref{intersect:h-intersect}.
The basic problem is to describe, classify, understand the strictly decreasing sequences 
$(A_q)_{q=1}^{\infty}$ of subsets of  
an additive abelian semigroup such that 
\beq           \label{intersect}
h\left( \bigcap_{q=1}^{\infty} A_q \right) = \bigcap_{q=1}^{\infty} hA_q.
\eeq
\eprob

For related work, see~\cite{fox-krav-zhan25}--\cite{schi25}.

\section{Intersections of sumsets in discrete groups and semigroups}

\bt               \label{intersect:theorem:finiteReps}
Let $h \in \N$.  Let $S$ be an additive abelian semigroup such that $r_{S,h}(x)$ 
is finite for all $x \in S$.  Let $A$ be a subset of $S$ 
and let $(A_q)_{q=1}^{\infty}$ be a decreasing sequence of subsets of $S$ 
such that $A = \bigcap_{q=1}^{\infty} A_q$. 
Then 
\[
hA = \bigcap_{q=1}^{\infty} hA_q.
\]
\et

\begin{proof}
For all $x \in S$, the set $ \mcr_{S,h}(x)$ of $h$-tuples $(b_1,\ldots, b_h) \in S^h$ 
such that 
\[
x = b_1 + \cdots + b_h 
\]
is finite with cardinality  $| \mcr_{S,h}(x)| = r_{S,h}(x)$.   

Let $x \in \bigcap_{q=1}^{\infty} hA_q$.  We must show that $x \in hA$. 
For all $q\in \N$, we have $x \in hA_q$ and so there is an $h$-tuple 
$(a_{1,q},\ldots, a_{h,q}) \in A_q^h \cap  \mcr_{S,h}(x)$ 
such that 
\[
x = a_{1,q} + \cdots + a_{h,q}. 
\]
Because $ \mcr_{S,h}(x)$ is a finite set, there is some $h$-tuple $(b_1,\ldots, b_h) \in  \mcr_{S,h}(x)$ 
such that 
\[
(b_1,\ldots, b_h) = (a_{1,q},\ldots, a_{h,q}) 
\]
for infinitely many $q$ and so, for infinitely many $q$, 
we have  $b_i = a_{i,q} \in A_q$ for all $i=1,\ldots, h$.  
Because the sequence $(A_q)_{q=1}^{\infty}$ is decreasing, 
we have $b_i \in A_q$ for all $i=1,\ldots, h$ and $q \in \N$.   
It follows that $b_i \in \bigcap_{q=1}^{\infty} A_q = A$ 
and $x = b_1 + \cdots + b_h \in hA$.
This completes the proof. 
\end{proof}

Let  $\Z^d$ be the additive abelian group of 
$d$-dimensional  lattice points in $\R^d$ 
and let $\N_0^d$ be the additive abelian semigroup of 
$d$-dimensional nonnegative lattice points. 

\bt               \label{intersect:theorem:nonnegative-Zd}
Let $d \geq 1$ and let $A$ be a  subset of $\N_0^d$. 
Let $(A_q)_{q=1}^{\infty}$ be a decreasing sequence of subsets of $\N_0^d$ 
such that $A = \bigcap_{q=1}^{\infty} A_q$. 
For every positive integer $h$, 
\[
hA = \bigcap_{q=1}^{\infty} hA_q.
\]
\et

\begin{proof} 
For all $n \in \N_0$, the number $r_{\N_0,h}(n)$  of representations 
of $n$ as a sum of $h$ nonnegative integers is finite.  
Indeed,  $r_{\N_0,h}(n) = \binom{n+h-1}{h-1}$.   
For all $\mbn = (n_1,\ldots, n_d) \in \N_0^d$, we have 
\[
r_{\N_0^d,h}(\mbn) = \prod_{i=1}^d r_{\N_o,h}(n_i) < \infty. 
\]
The result follows from Theorem~\ref{intersect:theorem:finiteReps}.  
\end{proof}

\bt                               \label{intersect:theorem:basis-S}
Let $S$ be an additive abelian semigroup and let $(A_q)_{q=1}^{\infty}$ 
be a decreasing sequence of subsets of $S$ such that 
$A = \bigcap_{q=1}^{\infty} A_q$.  
If $A$ is a basis of order $h_0$  for $S$, then 
\[
h_0A = \bigcap_{q=1}^{\infty} h_0 A_q 
\]
If the subgroup $S$ contains an additive identity and 
if $A$ is a basis of order $h_0$  for $S$, then 
\[
hA = \bigcap_{q=1}^{\infty} hA_q 
\]
 for all $h \geq h_0$. 
 \et
 
\begin{proof} 
For all $q \in \N$, we have $A \subseteq A_q$ and so $S = h_0A  \subseteq  h_0A_q  \subseteq  S$. 
It follows that $h_0A_q = S$ and 
\[
h_0A  = S =  \bigcap_{q=1}^{\infty} hA_q. 
\] 

Let $h \geq h_0$.  If $S$ contains an additive identity, then  $h_0A \subseteq  hA$. 
If $A$ is a basis of order $h_0$ for $S$, 
then  $S = h_0A \subseteq  hA \subseteq  S$ 
and so  $A$ is a basis of order $h$ for $S$.
This completes the proof. 
 \end{proof}

Note that the additive semigroup \N\ of positive integers does not contain an additive identity. 
The set $A= \N$ is a basis of order 1 for \N\ but not a basis of order $h$ for any $h \geq 2$.

\bt                   \label{intersect:theorem:Rd-nonbasis}
Let $d \in \N$ and let $G$ be any subgroup of  $\R^d$ that contains  $\Z^d$. 
Let $A$ be a subset of $G$.  
If $A$ is a nonbasis of order $h_0$  for $G$, then there exists a decreasing sequence 
$(A_q)_{q=1}^{\infty}$ of subsets of  $G$ 
such that $ A = \bigcap_{q=1}^{\infty} A_q$ and 
\[
 h_0A  \neq  \bigcap_{q=1}^{\infty} h_0 A_q.  
\] 
If $A$ is a nonbasis  for $G$, then there exists a decreasing sequence 
$(A_q)_{q=1}^{\infty}$ of subsets of  $G$ 
such that $ A = \bigcap_{q=1}^{\infty} A_q$ and 
\[
 hA  \neq  \bigcap_{q=1}^{\infty} hA_q  
\]
for all $h \geq 2$. 
\et

\begin{proof} 
For all $\mbx =  (x_1,\ldots, x_d) \in \R^d$, define $\| \mbx \| = \max( |x_i| : i = 1,\ldots, d)$.  
For all $q \in \N$, let 
\[
A_{q,d}^{\sharp} = \{ \mbx \in G: \| \mbx \| \geq q \}. 
\] 
Note that $\{ \mbx \in \Z^d: \| \mbx \| \geq q \} \subseteq A_{q,d}^{\sharp}$. 

Let $\mbx =  (x_1,\ldots, x_d) \in G$.  If $x_1 \geq 0$, then 
\[
\mby =  (x_1 + (h-1)q, x_2,\ldots, x_d) \in A_q^{\sharp} 
\]
and  
\[
\mbz = (-q,0,\ldots, 0) \in A_q^{\sharp}  
\]
and so 
\begin{align*}
\mbx & =  (x_1 + (h-1)q, x_2,\ldots, x_d) + (h-1)(-q,0,\ldots, 0) \\
& = \mby + (h-1)\mbz \in h A_{q,d}^{\sharp}.
\end{align*}
 If $x_1 \leq 0$, then 
 \[
 (x_1 - (h-1)q, x_2,\ldots, x_d) \in A_q^{\sharp} 
 \]
 and 
 \[
 (q,0,\ldots, 0) \in A_q^{\sharp}  
\]
and so 
\begin{align*}
\mbx & =  (x_1 - (h-1)q, x_2, \ldots, x_d) + (h-1)(q,0,\ldots, 0) \\
& = \mby + (h-1)\mbz \in h A_{q,d}^{\sharp}.
\end{align*}
Thus,  $hA_{q,d}^{\sharp} = G$ for all $q \in \N$ and $h \geq 2$. 

For all $q \in \N$, let 
\[ 
A_q = A \cup A_{q,d}^{\sharp}.   
\] 
Then  $hA_q = hA_{q,d}^{\sharp} = G$ for all $q \in \N$ and $h \geq 2$.   
The sequence  $(A_q)_{q=1}^{\infty}$ is a decreasing sequence of sets  
such that $ A = \bigcap_{q=1}^{\infty} A_q$.     
If $A$ is a nonbasis of order $h_0$  for $G$, then 
\[
h_0A \neq G = \bigcap_{q=1}^{\infty} h_0A_q.  
\]
 If $A$ is a nonbasis for $G$, then 
\[
 hA  \neq G =  \bigcap_{q=1}^{\infty} hA_q    
\] 
  for all $h \geq 2$.  This completes the proof. 
\end{proof}

\bt                 \label{intersect:theorem:Z-bounded below}
Let $A$ be an infinite set of integers that is bounded below and that does not contain 
all sufficiently large integers.  
There exists an asymptotically strictly decreasing 
sequence $(A_q)_{q=1}^{\infty}$ such that $hA = \bigcap_{q=1}^{\infty} hA_q$ for all $h \in \N$. 
\et 

\begin{proof}
Note that a set of integers that is bounded below cannot be a basis for \Z.  
By Theorem~\ref{intersect:theorem:nonnegative-Zd} with $d=1$, 
we can assume that $a^* = \min(A) < 0$.  

For all $q \in \N$, let 
 \[
 A^{\flat}_q =  \{r\in \N: r \geq q \} 
\] 
and let 
\[
A_q = A \cup  A^{\flat}_q = A \cup \{r\in \N: r \geq q\}.
\]
The sequence $(A_q)_{q=1}^{\infty}$ is decreasing and $A = \bigcap_{q=1}^{\infty} A_q$.  
If $q \in \N$ and $q \notin A$, then $q \in A_q$ but $q \notin A_{q+1}$ and so $A_q \neq A_{q+1}$.  
Because $A$ does not contain all sufficiently large integers, 
the sequence $(A_q)_{q=1}^{\infty}$ is asymptotically strictly decreasing.  

For all $h \in \N$ and $q \in \N$, we have 
\[
hA \subseteq hA_q \subseteq hA \cup \{r\in \N: r \geq q + (h-1)a^*\}
\]
and 
\[
\bigcap_{q=1}^{\infty}  \{r\in \N: r \geq q + (h-1)a^*\}  = \emptyset. 
\]
Therefore, 
\[
hA \subseteq  \bigcap_{q=1}^{\infty} hA_q
 \subseteq  hA \cup \bigcap_{q=1}^{\infty}  \{r\in \N: r \geq q + (h-1)a^*\} 
= hA.
\] 
This completes the proof. 
\end{proof}

Let $S$ be an additive abelian semigroup.  
A subset $A$ of $S$ is a \emph{maximal nonbasis of order $h$} 
for $S$ if  $A$  is not a basis of order $h$ for $S$ 
but every proper superset of $A$ is a basis of order $h$ for $S$.  
For example, for every integer $h \geq 2$,  the set 
$h\ast \Z = \{ h  \ell  : \ell \in \Z\}$ of multiples of $h$ is a maximal nonbasis 
of order $h$ for \Z\ if and only if $h$ is prime.    Maximal nonbases were introduced by 
Nathanson~\cite{nath74b} and the first nontrivial examples were constructed 
by Nathanson~\cite{nath74b} 
and Erd\H os and Nathanson~\cite{nath75a}.

\bt          \label{intersect:theorem:maximal}
Let $A$ be a subset of  the infinite abelian semigroup $S$.  
If $A$ is a maximal nonbasis of order $h$ 
for $S$ and if $(A_q)_{q=1}^{\infty}$ is an asymptotically  strictly decreasing sequence 
of subsets of $S$ such that $A = \bigcap_{q=1}^{\infty} A_q$, 
then $hA \neq  \bigcap_{q=1}^{\infty} hA_q$. 
\et

\begin{proof}
Let $A$ be a maximal nonbasis of order $h$ for $S$ and let  
$(A_q)_{q=1}^{\infty}$ be a  strictly decreasing sequence 
such that $A = \bigcap_{q=1}^{\infty} A_q$.  
For all $q$, the set $A$ is a proper subset of $A_q$ and so $hA_q = S$.  
It follows that $hA \neq S = \bigcap_{q=1}^{\infty} hA_q$.  
This completes the proof.  
\end{proof}

\bprob 
Let $A$ be a subset of  the infinite abelian semigroup $S$ 
such that $A$ is a  nonbasis of order $h$ but not a maximal nonbasis of order $h$ 
for $S$.   Does there exist an asymptotically  strictly decreasing sequence 
$(A_q)_{q=1}^{\infty}$ of subsets of $S$ such that $A = \bigcap_{q=1}^{\infty} A_q$
and $ hA = \bigcap_{q=1}^{\infty} h A_q$?  
In particular, consider this problem for $S = \Z$.  
\eprob

\section{Intersections of sumsets in locally compact groups} 

\bt                        \label{intersect:theorem:lca}
Let  $A$ be a nonempty  subset of a locally compact abelian group $G$.  
Let $(A_q)_{q=1}^{\infty}$ 
be a decreasing sequence of   compact subsets of $G$ 
such that $A = \bigcap_{q=1}^{\infty} A_q$. 
Then 
\[
hA = \bigcap_{q=1}^{\infty} hA_q 
\]
for every positive integer $h$.  
\et

\begin{proof}
For all $h \in \N$ and $q \in \N$, the sumset $hA_q$ and the sumset  intersection 
$\bigcap_{q=1}^{\infty} hA_q$ are nonempty closed compact sets.
Let $x \in \bigcap_{q=1}^{\infty} hA_q$.    We must show that $x \in hA$.  

For all $q\in \N$, we have $x \in hA_q$ and so there is an $h$-tuple  
\[
(a_{1,q},\ldots, a_{h,q}) \in A_q^h \subseteq A_1^h 
\]  
such that 
\[
x = a_{1,q} + \cdots + a_{h,q}. 
\]
Consider the sequences $(a_{i,q})_{q=1}^{\infty}$ for $i = 1, \ldots, h$. 
Because the set $A_1$ is compact and $a_{i,q} \in A_1$ for all  $i = 1, \ldots, h$ and $q\in \N$, 
there is an infinite strictly increasing sequence of positive integers $(q_j)_{j=1}^{\infty}$ such that 
the $h$ subsequences  $(a_{i,q_j})_{j=1}^{\infty}$ converge: 
\[
\lim_{j\rightarrow \infty} a_{i,q_j} = x_i   
\]
 for all  $i = 1, \ldots, h$.   Then 
\begin{align*}
x & = \lim_{j\rightarrow \infty} \left( a_{1,q_j} + \cdots + a_{h,q_j} \right) \\ 
& = \left( \lim_{j\rightarrow \infty}  a_{1,q_j}  \right)  + \cdots +  \left( \lim_{j\rightarrow \infty}  a_{h,q_j} \right) \\
& = x_1 + \cdots + x_h.
\end{align*}

Suppose that $x_i \notin A$ for some $i$.  Let $q_0=0$. 
Construct  the sequence $\left( x_{i,q} \right)_{q=1}^{\infty}$ as follows;
\[
x_{i,q} = a_{i,q_j} \qquad \text{for $q \in [q_{j-1}+1,q_j]$.} 
\]
Because $(A_q)_{q=1}^{\infty}$ is a decreasing sequence, 
we have $x_{i,q} \in A_q$ for all $q \in \N$ and 
\[
\lim_{i\rightarrow \infty} x_{i,q} = x_i.  
\]
The group $G$ is locally compact.  
Because $A = \bigcap_{q=1}^{\infty} A_q$ is closed and $x_i \notin A$, there is a compact neighborhood $K$  of $x_i$ 
that is disjoint from $A$ such that $x_{i,q}  \in K$ 
for all sufficiently large $q$.  
Then  $x_{i,q}  \in K \cap A_q$ for all sufficiently large $q$.  
It follows that, for all $Q \in \N$, there exists $q_1 \geq Q$ with 
$x_{i,q_1}  \in K \cap A_{q_1} \subseteq K \cap  A_q $ for all $q \in [1,Q]$, and so  
\[
x_{i,q_1}  \in \bigcap_{q=1}^Q \left( K \cap A_q \right) \neq \emptyset.  
\]
Because the sets $K \cap A_q$ are compact,  the finite intersection property implies 
\[
\emptyset \neq  \bigcap_{q=1}^{\infty}   \left( K \cap  A_q \right) 
=  K \cap  \left(  \bigcap_{q=1}^{\infty} A_q \right) 
= K \cap A = \emptyset
\]
which is absurd.  
Therefore, $x_i \in A$ for all $i  = 1, \ldots, h$ and $x \in hA$.  
This completes the proof. 
\end{proof}

\bt                       \label{intersect:theorem:lca-Rd}
Let $(A_q)_{q=1}^{\infty}$ be a decreasing sequence of nonempty compact sets 
 in $\R^d$ such that $A = \bigcap_{q=1}^{\infty} A_q \neq \emptyset$. 
For every positive integer $h$, 
\[
hA = \bigcap_{q=1}^{\infty} hA_q.
\]
\et

\begin{proof}
This follows from Theorem~\ref{intersect:theorem:lca} with $G = \R^d$.  
\end{proof}

\bt
Let $G$ be a locally compact abelian group with Haar measure $\mu$.
For $h \in \N$, let $\theta_h \geq 0$ and let $(\theta_{h,q})_{q=1}^{\infty}$ 
be a decreasing sequence of nonnegative real 
numbers such that $\lim_{q\rightarrow \infty} \theta_{h,q} = \theta_h$.  
Let $(A_q)_{q=1}^{\infty}$ be a decreasing sequence of compact subsets of $G$ 
such that $\mu(hA_q) = \theta_{h,q}$ for all $q \in \N$.  Let $A  = \bigcap_{q=1}^{\infty} A_q$.
Then $\mu(hA) = \theta_h$. 
\et

\begin{proof}
For all $h \in \N$, the sequence  $(hA_q)_{q=1}^{\infty}$ is  a decreasing sequence 
of compact subsets of $G$.  
By Theorem~\ref{intersect:theorem:lca}, 
\[
hA =   \bigcap_{q=1}^{\infty} hA_q 
\]
and so 
\[
\mu(hA) = \mu \left( \bigcap_{q=1}^{\infty} hA_q \right) = \lim_{q\rightarrow \infty} \mu(hA_q) 
 = \lim_{q\rightarrow \infty} \theta_{h,q}  = \theta_h. 
\]
This completes the proof.  
\end{proof}

\section{Additional problems}

\bprob
Let $(A_q)_{q=1}^{\infty}$ be a strictly decreasing sequence of sets 
in an additive abelian semigroup and let 
$A = \bigcap_{q=1}^{\infty}A_q$.  
Describe the set $\{h \in \N: hA = \bigcap_{q=1}^{\infty} hA_q\}$. 
What sets of integers arise in this way? 
\eprob

\bprob
Let $(A_q)_{q=1}^{\infty}$ be a strictly decreasing sequence of sets 
in an additive abelian semigroup and let 
$A = \bigcap_{q=1}^{\infty}A_q$.  \\
When does $hA = \bigcap_{q=1}^{\infty} hA_q$ 
imply $(h+1)A = \bigcap_{q=1}^{\infty} (h+1)A_q$?  \\
When does $hA = \bigcap_{q=1}^{\infty} hA_q$ 
imply $(h-1)A = \bigcap_{q=1}^{\infty} (h-1)A_q$?  
\eprob

\bprob
For what sets $A$ of integers does there exist a decreasing sequence 
$(A_q)_{q=1}^{\infty}$ of sets of  integers such that $A = \bigcap_{q=1}^{\infty} A_q$
and $hA \neq\bigcap_{q=1}^{\infty} hA_q$ for all $h \in \N$? 
\eprob

\emph{Acknowledgement.}  
I thank Diego Marques for helpful comments   that suggested  Theorem~\ref{intersect:theorem:maximal}.

\end{document}